# Summation arithmetic functions with asymptotically independent summands

VICTOR VOLFSON

ABSTRACT.

The summation arithmetic functions with asymptotically independent summands are studied in the paper. We prove statements about the condition under which the summation arithmetic functions have asymptotically independent summands. It is also prove that the limiting distribution of the summation arithmetic function with asymptotically independent summands is normal under certain conditions for summands of the summation arithmetic function.

1. INTRODUCTION

An arithmetic function (in the general case) is a function $f(n)$ defined on the set of natural numbers and taking values on the set of complex numbers. The name arithmetic function is related to the fact that this function expresses some arithmetic property of the natural series.

A summation arithmetic function is a function:

$$S(x) = \sum_{n \leq x} f(n). \tag{1.1}$$

Examples of summation functions (1.1) are: Mertens function $M(x) = \sum_{k \leq x} \mu(k)$, where $\mu(k)$ is Möbius function, the summation function of the number of square-free numbers:

$$Q(x) = \sum_{k \leq x} |\mu(k)| \tag{1.2}$$

and others.

A Mobius function $\mu(k) = 1$ if the natural number k has an even number of prime divisors of the first degree, $\mu(k) = -1$ if the natural number k has an odd number of prime divisors of the first degree and $\mu(k) = 0$ if the natural number k has prime divisors of not only first degree.

______________________________________________________________________

The sequence of the values of the summation arithmetic function $S(1),...,S(n)$ can be distributed altogether quite oddly. If you follow the values of such functions when the argument runs through the natural values, you get a very chaotic picture.

In classical studies, the average value of arithmetic functions is usually considered when studying their distributions arithmetic on the initial segment of the natural series $1,...,n$ and find asymptotic approximate expressions for it [1].

If we want to more accurately characterize the distribution of values $S(n)$ on the initial segment of the natural series $1,...,n$, then we come to the concept of the distribution function $S(n)$ - $P(S(n) < x) = F_n(x)$.

Since $F_n(x)$ is a distribution function $S(n)$ on the initial segment of the natural series $1,...,n$ in a probability-theoretic sense, it is natural to consider the convergence of the sequence of distribution functions $F_n(x)$ for $S(n)$ on the entire natural series to the limit distribution function under certain conditions.

Definition. We will understand as asymptotic independence of arithmetic functions, which the limit of the difference between the average value of the product of an arithmetic function $f(k)$ and the product of the average values of the same function at different values of the argument tends to zero when $n \to \infty$.

It was shown in [2] using probabilistic methods, that the asymptotic independence of the summands in combination with other properties leads to a normal limiting distribution of the summation function $S(x) = \sum_{n \leq x} f(n)$.

It was possible to prove in [3], [4] already with the help of exact methods, that the property of asymptotic independence of summands is fulfilled for the summation functions of Mertens $M(x) = \sum_{n \leq x} \mu(n)$ and Liouville $L(x) = \sum_{n \leq x} \lambda(n)$. We indicate the property of the summation arithmetic functions in the paper, which is sufficient to fulfill the asymptotic independence of their terms. We consider summation arithmetic functions for which the asymptotic independence of the summands will be satisfied, including the summation function of the number of square-free numbers.

In the paper, we will significantly expand the class of summation arithmetic functions for which the asymptotic independence of the summands of arithmetic functions holds.

Further, we show that summation arithmetic functions which are composed bounded arithmetic functions (for which asymptotic independence is fulfilled) have a limit normal distribution under certain additional conditions.

2. ASYMPTOTIC INDEPENDENSE OF THE SUMMATION ARITHMETIC FUNCTIONS

Denote Möbius or Liouville arithmetic functions - $f(k)$.

Statement 1

Let the average value of the product of the arithmetic function with different values of the argument is determined by the formula:

$$\frac{\sum_{i=1}^{n} \sum_{j=1(i \neq j)}^{n} f(i)f(j)}{n(n-1)} . \qquad (2.1)$$

Let the product of the average values of the arithmetic function for different values of the argument is determined by the formula:

$$\frac{(\sum_{k=1}^{n} f(k))^2 - \sum_{k=1}^{n} f^2(k)}{n^2} . \qquad (2.2)$$

Then, the upper bound for the difference between the average value of the product of an arithmetic function $f(k)$ and the product of the average values of the same function for different values of the argument is $o(1/n)$.

Proof

Find the difference:

$$\frac{\sum_{i=1}^{n}\sum_{j=1(i\neq j)}^{n} f(i)f(j)}{n(n-1)} - \frac{(\sum_{k=1}^{n} f(k))^2 - \sum_{k=1}^{n} f^2(k)}{n^2}. \tag{2.3}$$

Having in mind $\sum_{i=1}^{n}\sum_{j=1(i\neq j)}^{n} f(i)f(j) = \sum_{i=1}^{n}\sum_{j=1}^{n} f(i)f(j) - \sum_{k=1}^{n} f^2(k)$ and substituting it in (2.3), we get:

$$\frac{\sum_{i=1}^{n}\sum_{j=1}^{n} f(i)f(j) - \sum_{k=1}^{n} f^2(k)}{n(n-1)} - \frac{(\sum_{k=1}^{n} f(k))^2 - \sum_{k=1}^{n} f^2(k)}{n^2}. \tag{2.4}$$

Based on (2.4), having in mind $\sum_{i=1}^{n}\sum_{j=1}^{n} f(i)f(j) = (\sum_{k=1}^{n} f(k))^2$, we get the estimate:

$$(\sum_{k=1}^{n} f(k))^2 - \sum_{k=1}^{n} f^2(k))(1/n(n-1) - 1/n^2). \tag{2.5}$$

Having in mind $(\sum_{k=1}^{n} f(k))^2 = o(n^2)$ and $\sum_{k=1}^{n} f^2(k) = O(n)$, substituting it in (2.5), we get the desired estimate: $o(1/n)$

Corollary 1

The asymptotic independence for Möbius or Liouville arithmetic functions is observed.

Proof

Based on Theorem 1, the limit of difference between the average value of the product of an arithmetic function $f(k)$ and the product of the average values of the same function at different values of the argument tends to zero when $n \to \infty$, therefore, the asymptotic independence of Mobius or Liouville arithmetic functions is fulfilled.

An indicator of asymptotic independence of the summands of arithmetic functions is the asymptotic upper bound for the difference of expressions (2.1) and (2.2). For example, asymptotic independence is performed with the exponent $o(1/n)$ for the functions indicated in Corollary 1.

The summing function — the number of primes not exceeding $n$ - $\pi(n)$ has asymptotic $\pi(n) \sim \dfrac{n}{\log(n)} = o(n)$, therefore, for its terms, the relation is satisfied - $(\sum_{k=1}^{n} f(k))^2 = o(n^2)$.

On the other hand $\pi(n) = \sum_{p} 1$, therefore, the terms of this summation arithmetic function are bounded and the relation holds - $\sum_{k=1}^{n} f^2(k) = O(n)$.

Consequently, based on (2.5), the asymptotic independence of the terms for the summation arithmetic function $\pi(n)$ is also performed with the exponent $o(1/n)$.

Corollary 2

Asymptotic independence is performed with the exponent $O(1/n)$ for bounded arithmetic functions - $f(n)$.

Proof

It is performed $(\sum_{k=1}^{n} f(k))^2 = O(n^2), \sum_{k=1}^{n} f^2(k) = O(n)$ for bounded arithmetic functions; therefore, having in mind (2.5), asymptotic independence is performed with the exponent - $O(1/n)$.

However, the property of asymptotic independence holds also for summation arithmetic functions whose summands are not limited.

Statement 2

If the summation arithmetic function $S(x) = \sum_{k=1}^{n} f(k)$ has the summands of the function $f(k)$ of the same sign, then it is executed:

$$(\sum_{k=1}^{n} f(k))^2 - \sum_{k=1}^{n} f^2(k))(1/n(n-1) - 1/n^2) = \dfrac{O(S^2(n))}{n^3}. \tag{2.6}$$

Proof

Based on (2.5):

$$(\sum_{k=1}^{n} f(k))^2 - \sum_{k=1}^{n} f^2(k))(1/n(n-1) - 1/n^2) = \frac{O((\sum_{k=1}^{n} f(k))^2 + O(\sum_{k=1}^{n} f^2(k)))}{n^3}. \qquad (2.7)$$

Having in mind:

$$S^2(n) = (\sum_{k=1}^{n} f(k))^2 = \sum_{k=1}^{n} f^2(k) + \sum_{i=1}^{n} \sum_{j=1,(i \neq j)}^{n} f(i)f(j). \qquad (2.8)$$

Based on (2.8), if the terms $f(k)$ are of the same sign, then:

$$\sum_{k=1}^{n} f^2(k) = O(S^2(n)). \qquad (2.9)$$

Having in mind (2.7) in (2.9) we get:

$$\sum_{k=1}^{n} f^2(k))(1/n(n-1) - 1/n^2) = \frac{O(S^2(n))}{n^3} + \frac{O(S^2(n))}{n^3} = \frac{O(S^2(n))}{n^3},$$

which corresponds to (2.6).

Note, that the proof did not use no conditions of the limitations $f(k)$.

Corollary 3

If $S(n) = O(n)$ and the summands of the summation arithmetic function do not change sign, then based on (2.6):

$$(\sum_{k=1}^{n} f(k))^2 - \sum_{k=1}^{n} f^2(k))(1/n(n-1) - 1/n^2) = O(1/n). \qquad (2.10)$$

The following summation arithmetic functions satisfy condition (2.10): Chebyshev - $\Psi(x) = x + o(x), \theta(x) = x + o(x)$, $R(x) = \sum_{n \leq x} r(n) = \frac{\zeta(2)\zeta(3)}{\zeta(6)} x + o(x)$ (where $r(n)$ is the number of solutions of the equation $\varphi(m) = n$ and $\varphi(m)$ is the Euler function), the number of square-free numbers - $Q(x) = 6x/\pi^2 + o(x)$, the number of numbers having an odd number of prime dividers - $Q_1(x) = 3x/\pi^2 + o(x)$, the number of numbers having an even number of prime dividers - $Q_2(x) = 3x/\pi^2 + o(x)$.

Corollary 4

Asymptotic independence (with exponent $o(1)$) is performed for summation arithmetic functions satisfying the estimate - $S(x) = o(x^{3/2})$.

We must substitute the estimate $S(x) = o(x^{3/2})$ in formula (2.6) for proof.

Corollary 4 satisfies for summation function of the number of divisors of natural numbers not exceeding $x$ - $S(x) = \sum_{n \leq x} \tau(n) = x \log x + (2C-1)x + o(x)$. Therefore, the asymptotic independence of the arithmetic function $\tau(n)$ is performed (the number of positive divisors of the natural - $n$).

It can be shown that Corollary 4 also satisfies alternating arithmetic functions $f(n)$.

Thus, asymptotic independence of the summands of arithmetic functions is performed, respectively, with the indicators: $o(1/x), O(1/x), o(1)$, if the asymptotic upper bound for summation arithmetic functions is respectively: $o(x), O(x), o(x^{3/2})$.

3. LIMIT DISTRIBUTION OF SUMMATION ARITHMETIC FUNCTION

It is known [5] that any initial segment of a natural series $\{1, 2, ..., n\}$ can be naturally transformed into a probabilistic space $(\Omega_n, \mathcal{A}_n, \mathbb{P}_n)$, by taking $\Omega_n = \{1, 2, ..., n\}$, $\mathcal{A}_n$ - all subsets $\Omega_n$, $P_n(A) = \frac{1}{n}\{N(m \in A)\}$, where $N(m \in A)$ is the number of natural series members that satisfy the condition $m \in A$. Then an arbitrary (real) function of the natural argument $f(k)$ (or rather, its restriction to $\Omega_n$) can be considered as a random variable $\xi_n$ on this probability space: $\xi_n(k) = f(k)$. In particular, we can talk about the average - $M[\xi_n] = \frac{1}{n}\sum_{k=1}^{n} f(k)$, variance $D[\xi_n] = \frac{1}{n}\sum_{k=1}^{n}|f(k)|^2 - |\frac{1}{n}\sum_{k=1}^{n}f(k)|^2$, distribution function $F_{\xi_n}(x) = \frac{1}{n}\{k \leqslant n : f(k) \leqslant x\}$ and characteristic function $\varphi_{\xi_n}(t) = \frac{1}{n}\sum_{k=1}^{n}e^{itf(k)}$.

We consider random variables $f_k$ that (in accordance with the definition of probability spaces) take values equal to the values of the arithmetic function, therefore $f_k = f(k), (1 \leq k \leq n)$.

We assume that random variables $f_k$ are quasi-asymptotically independent if the following relationship holds for their average values when $n \to \infty$:

$$M_{ij}[f,n] \to M_i[f,n]M_j[f,n],$$

where $M_{ij}[f,n]$ is determined by formula (2.1), and $M_i[f,n]M_j[f,n]$ is determined by formula (2.2). Quasi, since random variables $f_k$ are defined in different probability spaces.

Let us consider the summation arithmetic function $S(n) = \sum_{k=1}^{n} f(k)$ and define a random variable - $S_n(k) = S(k), (1 \leq k \leq n)$.

We will look for the limiting distribution function for a random variable $S_n$ when the value $n \to \infty$. In accordance with the above definition of the probability space for any arithmetic function, we can say that our goal is to find the limit distribution function for a sequence of corresponding random variables $S_n$ when $n \to \infty$.

Naturally, not every random variable $S_n$ has a limit distribution function when $n \to \infty$. It is necessary for it that the conditions of the theorem on the continuity of the characteristic function be satisfied. We will look for the characteristic function for this random variable - $\varphi_{S_n}(t) = M[e^{itS_n}]$, provided that this theorem is fulfilled for $S_n$.

However, we cannot directly use the quasi-asymptotic independence of random variables $f_k, (1 \leq k \leq n)$ to find the characteristic function $S_n$, since random variables $f_1,...,f_n$ are defined in different probability spaces. This problem can be solved by the following statement.

Statement 3

If the random variables $f_k = f(k), (1 \leq k \leq n)$ are quasi-asymptotically independent when $n \to \infty$, and for each value $t$ there is a limit of the characteristic function - $lim_{n\to\infty}\varphi_{S_n}(t)$ continuous at the point $t=0$, then the following is true:

$$\varphi_{S_n}(t) = \prod_{k=1}^{n} \varphi_{f_k}(t) \qquad (3.1)$$

when $n \to \infty$ and $\varphi_{S_n}(t)$ uniquely determines the limiting distribution function for $S_n$.

Proof

Having in mind that random variables $f_k$ are quasi-asymptotically independent, then for their average values the above relation is fulfilled when $n \to \infty$:

$$M_{ij}[f,n] \to M_i[f,n]M_j[f,n], \qquad (3.2)$$

where $M_{ij}[f,n]$ is determined by formula (2.1) and $M_i[f,n]M_j[f,n]$ is determined by formula (2.2).

Based on Lemma 3, p. 123 [6], it is possible to construct random variables: $g_1,...,g_n,G_n$, which have, respectively, equal distribution functions with $f_1,...,f_n,S_n$, and hence characteristic functions. Since the distribution functions coincide, then also coincide the average values, respectively, for $f_k$ and $g_k$; therefore, we can write a similar relation for the average values for $g_k$ when performed (3.2) and when $n \to \infty$:

$$M_{ij}[g,n] \to M_i[g,n]M_j[g,n]. \qquad (3.3)$$

Based on (3.3), without taking into account the trivial cases, we can assume that $g_1,...,g_n$ are already asymptotically independent (without quasi), since they are in the same probability space. Based on this and the properties of the characteristic function we obtain (when $n \to \infty$):

$$\varphi_{\sum_{k=1}^{n} g_k}(t) = \prod_{k=1}^{n} \varphi_{g_k}(t) = \prod_{k=1}^{n} \varphi_{f_k}(t), \qquad (3.4)$$

having in mind the equality of the corresponding characteristic functions $\varphi_{g_k}(t) = \varphi_{f_k}(t)$.

Having in mind (3.4) and the condition that for each value there is a limit of the characteristic function - $lim_{n\to\infty}\varphi_{S_n}(t)$, continuous at the point $t=0$, and also based on the theorem on the continuity of the characteristic function, the product $\prod_{k=1}^{n}\varphi_{f_k}(t)$ is a characteristic function of a certain probability distribution for $\sum_{k=1}^{n}g_k$ when $n\to\infty$.

Let this limit probability distribution have a distribution function (limit distribution function) - $G(x)$.

Having in mind that the distribution functions for $f_k, g_k$, respectively, coincide, the distribution functions for $G_n, S_n$ and, respectively, the limit distribution functions also coincide.

Consequently, a random variable $S_n$ has a limit distribution function $G(x)$ when $n\to\infty$, which is uniquely determined by the characteristic function $\varphi_{\sum_{k=1}^{n}g_k}(t)$, which in turn, based on (3.4), is equal to $\prod_{k=1}^{n}\varphi_{f_k}(t)$.

Thus, when $n\to\infty$ it is satisfied $\varphi_{S_n}(t)=\prod_{k=1}^{n}\varphi_{f_k}(t)$, which corresponds to (3.1).

Having in mind that by the condition at each value $t$ there is a limit of the characteristic function - $lim_{n\to\infty}\varphi_{S_n}(t)$, continuous at a point $t=0$, based on the continuity theorem for the characteristic function $\varphi_{S_n}(t)$, the limiting distribution function for a random variable $S_n$ is uniquely determined. The proof is complete.

The limit distribution function for Möbius function was found in [7], and the limit distribution function for Liouville function was found in [4].

Now consider the more general case.

Statement 4

Let the arithmetic function $f: N \to R$ takes the values: $a_1, ..., a_k$. We define a random variable $f_n: f_n(m) = f(m), (1 \leq m \leq n)$ in a probability space - $(Q_n, A_n, P_n)$, where $Q_n=\{1,...,n\}$, $A_n$ are all the subsets $Q_n$ and $P_n : A_n \to R$. Let consider the values of probabilities: $v_1(n) = N(f(i) = a_1)/n$, ..., $v_k(n) = N(f(i) = a_k)/n$ where $1 \leq i \leq n$ and $v_1(n) + ... + v_k(n) = 1$.

Then, if, $\lim_{n \to \infty} v_1(n) = p_1, ..., \lim_{n \to \infty} v_k(n) = p_k$ then:

$$\lim_{n \to \infty} P(f_n < y) = G(y), \tag{3.5}$$

where $G_n(y) = P(f_n < y)$ is the distribution function of a random variable $f_n$, and $G(y)$ is the limiting distribution function for $G_n(y)$, which is equal to:

$$G(y) = \{0, y < a_1; p_1, a_1 \leq y < a_2; ...; p_1 + ... + p_{k-1}, a_{k-1} \leq y < a_k; 1, y \geq 1\}. \tag{3.6}$$

Proof

Based on the definition of a random variable $f_n$, its distribution function is equal to:

$$G_n(y) = \{0, y < a_1; v_1(n), a_1 \leq y < a_2; ...; v_1(n) + ... + v_{k-1}, a_{k-1} \leq y < a_k; 1, y \geq 1\}. \tag{3.7}$$

Having in mind (3.7) and Remarks 4 on p. 123 [6], the distribution functions $G_n(y)$ converge to the distribution function $G(y)$, when value $n \to \infty$, as having jumps in the same points.

Consequently, $\lim_{n \to \infty} P(f_n < y) = G(y)$, that corresponds to (3.5) and (3.6).

Corollary

A similar statement holds if the arithmetic function takes an infinite (countable) number of values.

Denote the random variable with the limit distribution function $G(y)$ - $f$. This designation will be used below.

The following arithmetic functions with asymptotically independent terms satisfy the conditions of Statement 4: Liouville function with the limit distribution function - $G(y) = \{0, y < -1; 0, 5, -1 \leq y < 1; 1, y \geq 1\}$, Möbius function with the limit distribution function - $G(y) = \{0, y < -1; 3/\pi^2, -1 \leq y \leq 0; 1 - 3/\pi^2, 0 \leq y < 1; 1, y \geq 1\}$ and others.

The terms for the summation arithmetic function - the number of primes not exceeding

$n$ - $\pi(n) = \sum_{k=1, k \in p}^{n} 1$ have a degenerate distribution function - $G(y) = \{0, y < 0; 1, y \geq 0\}$.

Statement 5

Let the random variable $f$ is bounded. Then, in the neighborhood, $t = 0$ the following Taylor expansion for the characteristic function $f$ holds true:

$$\varphi_f(t) = 1 + \sum_{j=1}^{l} \frac{(it)^j}{j!} M[f^j] + o(t^l), \qquad (3.8)$$

where $M[f^j]$ is the $j$ - th moment from the arithmetic function $f$.

Proof

The average value (mathematical expectation) of a random variable is limited - $M[f] < \infty$, variance, and other moments of higher orders of a random variable $f$ are also limited - $M[f^l] < \infty$ for a limited random variable.

Therefore, based on property 5 of the characteristic function on page 131 [6], the following is performed:

$$\varphi_f^{(j)}(0) = t^{j-1} M[f^{j-1}]. \qquad (3.9)$$

Consequently, since $M[f^l] < \infty$, having in mind (3.9), then in the neighborhood $t = 0$, the decomposition of the characteristic function of $f$ in a Taylor series is true:

$$\varphi_f(t) = 1 + \sum_{j=1}^{l} \frac{(it)^j}{j!} M[f^j] + o(t^l)$$,

which corresponds to (3.8).

This is a purely probabilistic statement. Its meaning is that if there is a bounded random variable, then its characteristic function has the indicated expansion in a Taylor series in a neighborhood of a point $t = 0$.

Statement 6

Suppose there are random bounded quantities $f_n$, and a random quantity $f : f_n \to f$ (by distribution) when the value $n \to \infty$.

Suppose that it is performed for the mathematical expectation $f_n$:

$$M[f_n] = M[f] + o(g(n)), \quad (3.10)$$

where $g(n)$ is a decreasing function and $lim_{n \to \infty} g(n) = 0$.

Then the asymptotic of the characteristic function $f_n$ in a neighborhood $t = 0$ when $n \to \infty$ is:

$$\varphi_{f_n}(t) = \varphi_f(t) + r, \quad (3.11)$$

where $|r| = |t| o(g(n))$.

Proof

Having in mind that the random variables $f_n$ are limited, the random variable $f : f_n \to f$ (by distribution) with the value $n \to \infty$ is also limited. Therefore, based on Statement 5, we obtain:

$$|r| = |\varphi_{f_n}(t) - \varphi_f(t)| = |it(M[f_n] - M[f]) - t^2/2(M[f_n^2] - M[f^2]) + ...| = |t| o(g(n)),$$

that corresponds to (3.11).

Statement 7

Suppose there is a summation function $S(n) = \sum_{k=1}^{n} f(k)$, where the term arithmetic function $f : N \to R$ is bounded. Suppose that all the conditions of Statement 6 are satisfied for a random variable $f_n : f_n(k) = f(k), (n = 1, 2, ...)$ with $g = 1/n$. Then the random variable $S_n : S_n(k) = S(k), (1 \le k \le n)$, and therefore the summation function $S(n)$ has a normal distribution when $n \to \infty$.

Proof

The arithmetic function $f: N \to R$ is bounded by the condition, and therefore is asymptotically independent (Statement 1, Corollary 2). Based on Statement 3, random variables $f_k, (1 \leq k \leq n)$ are quasi-asymptotically independent if, for each value $t$, there is a limit of the characteristic function $lim_{n \to \infty} \varphi_{S_n}(t)$ that is continuous at the point $t = 0$ (we will show this later).

Therefore, it is executed when $n \to \infty$:

$$\varphi_{S_n}(t) = \prod_{k=1}^{n} \varphi_{f_k}(t). \tag{3.12}$$

Based on Statement 6, we have:

$$\varphi_{S_n}(t) = (\varphi_f(t) + r)^n = (\varphi_f(t))^n + R, \tag{3.13}$$

where $R = nr\varphi_f(t)^{n-1}) + n(n-1)/2(r^2\varphi_f(t)^{n-2)}) + ... + r^n$.

Based on (3.13) having in mind that for characteristic function holds $|\varphi_f(t)| \leq 1$ we get:

$$|R| \leq |nr| + |r^2 n(n-1)/2| + ... + |r|^n. \tag{3.14}$$

Based on Statement 6 and (3.14), we get when value $n \to \infty$:

$$|R| = no(1/n)|t| = o(1)|t|. \tag{3.15}$$

Having in mind (3.15) the value $R \to 0$ when $n \to \infty$.

Let us denote the mathematical expectation of a random variable $f$ - $m$, and the variance - $\sigma^2$, then based on Statement 5 (by assumption the random variable is limited) we get:

$$\varphi_f(t) = 1 + itm - \sigma^2 t^2/2 + o(t^2). \tag{3.16}$$

Instead of a random variable $S_n$, we consider a random variable $Z_n = (S_n - mn)/\sigma\sqrt{n}$. Then, to prove the statement, it is necessary to show that the value tends to $\varphi_{Z_n}(t) \to e^{-t^2/2}$ when $n \to \infty$ and $m = 0$.

Based on (3.12) - (3.16) when $m = 0$ we get:

$$\varphi_{Z_n}(t) = \varphi_f^n(t/\sigma\sqrt{n}) + R = [1 - t^2/2n + o(t^2/\sigma^2 n)]^n + R. \qquad (3.17)$$

Having in mind (3.17) we get:

$$\lim_{n\to\infty}\log[\varphi_{Z_n}(t) - R] = \lim_{n\to\infty}\log[\varphi_{Z_n}(t)] = -t^2/2. \qquad (3.18)$$

Therefore, based on (3.18) when $n \to \infty$ we get:

$$\varphi_{Z_n}(t) \to e^{-t^2/2}.$$

The function $e^{-t^2/2}$ is a characteristic function of the normal distribution, which, by virtue of the theorem on the continuity of the characteristic function, proves this statement.

Let us explain the above with examples.

Example 1

We consider the summation functions of Mertens and Liouville. Let us use the notation adopted above for the indicated summation functions: $S(n) = \sum_{k=1}^{n} f(k)$.

Let us check the conditions of Statement 7 for these functions.

The terms of the function $f: N \to R$ (Mobius or Liouville) and, accordingly, the random variables $f_n: f_n(k) = f(k)$ and the random variable $f: f_n \to f$ (based on Statement 4) are bounded; therefore, this part of the condition of Statement 7 is satisfied.

The following estimate is performed for these functions - $S(n) = o(n^{1/2+\epsilon})$ (where $\epsilon$ is a small positive number), if the Riemann hypothesis is true. Therefore, the estimate for the difference between the average values of the terms of arithmetic functions and, accordingly, for random variables, is:

$$M[f_n] - M[f] = 1/n \sum_{k=1}^{n} f(k) = S(n)/n = o(\frac{1}{n^{1/2+\epsilon}}), \qquad (3.19)$$

having in mind that the value $M[f] = 0$.

The function $1/n$ has an order of smallness higher than $M[f_n]$ in (3.19), therefore this part of Statement 7 is not satisfied and we cannot say that the summation arithmetic functions of Mertens or Liouville have a limit normal distribution when $n \to \infty$.

We formulate Example 2 as a statement.

Statement 8

Let us consider the alternating convergent series, in which both the series of positive terms and negative terms converge. The sum of the series of positive terms is equal to $\sum_{k=1}^{\infty} a_k = A$, and the sum of the series of negative terms is equal to $\sum_{k=1}^{\infty} b_k = -A$. Then the summation function, which is a partial sum of this alternating series $S(n) = \sum_{k=1}^{n} f(k) = \sum_{k=1}^{n} (a_k + b_k)$, has a limit normal distribution when $n \to \infty$.

Proof

Having in mind the condition the sum of the alternating series is equal to:

$$\sum_{k=1}^{\infty} (a_k + b_k) = A - A = 0. \qquad (3.20)$$

Since the series $\sum_{k=1}^{\infty} a_k$ converges, then $a_n \to 0$ when $n \to \infty$. Since the series $\sum_{k=1}^{\infty} b_k$ converges, then $b_n \to 0$ when $n \to \infty$. Therefore:

$$a_n + b_n \to 0 \text{ when } n \to \infty \qquad (3.21)$$

It is follows the boundedness of the arithmetic $f(n) = a_n + b_n$ based on (3.21), therefore, the first part of the condition of Statement 7 is satisfied.

Having in mind (3.21) the random variable $f = 0, (f_n \to f)$ (where the random variable $f_n : f_n(k) = f(k)$), therefore, $f$ is limited and $M[f] = 0$.

Since the series of positive terms converge, then the estimate is valid for a partial amount $\sum_{k=1}^{n} a_k = A + o(1)$. Since the series of negative terms converges, then the estimate is valid for a partial amount $\sum_{k=1}^{n} b_k = -A + o(1)$. Therefore:

$$\sum_{k=1}^{n}(a_k + b_k) = A + o(1) - A + o(1) = o(1). \quad (3.22)$$

Having in mind (3.22), we obtain the following estimate for the summation arithmetic function:

$$S(n) = \sum_{k=1}^{n} f(k) = \sum_{k=1}^{n}(a_k + b_k) = o(1). \quad (3.23)$$

Based on (3.21) and (3.23) we obtain the estimate for the average value of a random variable:

$$M[f_n] = \sum_{k=1}^{n} f(k)/n = M[f] + o(1/n), \quad (3.24)$$

which corresponds to the second part of Statement 7.

All the conditions of Statement 7 are fulfilled, therefore this summation arithmetic function has a limit normal distribution when $n \to \infty$.

4. CONCLUSION AND SUGGESTIONS FOR FURTHER WORK

The next article will continue to study the behavior of arithmetic functions.

5. ACKNOWLEDGEMENTS

I am grateful to the dxdy scientific forum (hhttps://dxdy.ru/) for valuable comments and suggestions made at the discussion of this work, which have significantly improved its content.

References

1. A.A. Buchstab. "Number Theory", Izvestia "Enlightenment", Moscow, 1966, 384 p.

2. Volfson V. Investigations of the limit distribution and the asymptotic behavior of summation arithmetic functions, arXiv preprint https://arxiv.org/abs/1804.07539 (2018)

3. Volfson V. Comparison of probabilistic and exact methods for estimating the asymptotic behavior of summation arithmetic functions, arXiv preprint https://arxiv.org/abs/1805.11880 (2018)

4. Volfson V. Investigations of Mertens and Liouville summation functions, arXiv preprint https://arxiv.org/abs/1807.09486 (2018)

5. J. Kubilius - Probabilistic methods in number theory, 220 p., Vilnius, 1962.

6. Borovkov A.A. Probability theory. M .: Editorial URSS, 1999,472 p.

7. Volfson V. Investigation of the asymptotic behavior of the Mertens function, arXiv preprint https://arxiv.org/abs/1712.04674 (2017)

8. Shiryaev A.N. Probability is 2. Moscow: Publishing house "Enlightenment 1966,384 p.